\documentclass[11pt]{article}
\usepackage{epic,latexsym,amssymb}
\usepackage{amsfonts}
\usepackage{amscd}
\usepackage{amsmath}
\usepackage{graphicx}
\usepackage{color}
\usepackage{caption,
caption}
\usepackage{tikz}
\usetikzlibrary{arrows.meta}
\usepackage{algorithm}
\usepackage[algo2e,ruled,vlined]{algorithm2e}

\usepackage{float}

\usepackage{graphicx}

\usepackage[utf8]{inputenc}
\usepackage{longtable} 

\newtheorem{duplicate}{Theorem}

\newtheorem{thm}{Theorem}
\newtheorem{conj}{Conjecture}
\newtheorem{ob}[thm]{Observation}

\newtheorem{lem}[thm]{Lemma}
\newtheorem{quest}{Question}

\newtheorem{cor}[thm]{Corollary}
\newtheorem{prob}{Problem}
\newtheorem{example}{Example}
\newtheorem{remark}{Remark}

\newcommand{\core}{{\rm core}}

\newcommand{\proof}{\noindent\textbf{Proof. }}

\newcommand{\qed}{$\Box$}

\newcommand{\1}{\vspace{0.1cm}}

\newcommand{\QEDmark}{\mbox{\textsc{qed}}}
\newcommand{\proofStarter}[1]{\textsc{#1} }

\def\vertex(#1){\put(#1){\circle*{2}}}
\def\vertexo(#1){\put(#1){\circle{2}}}
\def\vert(#1){\put(#1){\circle*{1.5}}}
\def\verto(#1){\put(#1){\circle{1.5}}}
\def\lab(#1)#2{\put(#1){\makebox(0,0)[c]{#2}}}
\definecolor{DarkGreen}{rgb}{0.2, 0.6, 0.3}

\definecolor{electricindigo}{rgb}{0.44, 0.0, 1.0}

\let\oldenumerate\enumerate
\renewcommand{\enumerate}{
  \oldenumerate
  \setlength{\itemsep}{0.5pt}
  \setlength{\parskip}{0pt}
  \setlength{\parsep}{0pt}
}

\makeatletter
\newcommand{\Rmnum}[1]{\expandafter\@slowromancap\romannumeral #1@}
\makeatother

\begin{document}

\title{New results relating independence \\ and matchings}

\author{$^{1}$Yair Caro, $^{2}$Randy Davila, and $^{2}$Ryan Pepper \\
\\
$^1$Department of Mathematics\\
University of Haifa--Oranim\\
Tivon 36006, Israel \\
\small {\tt Email: yacaro@kvgeva.org.il}
\\
\\
$^2$Department of Mathematics and Statistics\\
University of Houston--Downtown \\
Houston, TX 77002 \\
\small {\tt Email: davilar@uhd.edu} \\
\small {\tt Email: pepperr@uhd.edu} 
}

\date{}
\maketitle

\begin{abstract}
In this paper we study relationships between the \emph{matching number}, written $\mu(G)$, and the \emph{independence number}, written $\alpha(G)$. Our first main result is to show
\[
\alpha(G) \le \mu(G) +  |X| - \mu(G[N_G[X]]),
\]
where $X$ is \emph{any} intersection of maximum independent sets in $G$. Our second main result is to show
\[
\delta(G)\alpha(G) \le \Delta(G)\mu(G),
\] 
where $\delta(G)$ and $\Delta(G)$ denote the minimum and maximum vertex degrees of $G$, respectively. These results improve on and generalize known relations between $\mu(G)$ and $\alpha(G)$. Further, we also give examples showing these improvements. 
\end{abstract}

\date{}
\maketitle


\section{Introduction}
Graphs considered here will be finite, undirected, and with no loops. Let $G$ be a graph with order $n(G) = |V(G)|$ and size $m(G) = |E(G)|$. The open neighborhood of a vertex $v\in V(G)$ is the set of all vertices adjacent to $v$, written $N_G(v)$, whereas the closed neighborhood of $v$ is $N_G[v] = N_G(v) \cup \{v\}$. The minimum and maximum vertex degrees of $G$ will be denoted $\delta(G)$ and $\Delta(G)$, respectively. 

For a subset $X \subseteq V(G)$, we will use the notations $N_G(X) = \bigcup_{v\in X}N_G(v)$ and $N_G[X] = X \cup N_G(X)$, also $G[X]$ will denote the subgraph induced by $X$. A matching is a subset $M\subseteq E(G)$ of non-adjacent edges. Vertices incident with a matching are called \emph{saturated} by that matching. The \emph{matching number} is the cardinality of a maximum matching in $G$, and will be denoted by $\mu(G)$. A subset $X\subseteq V(G)$ is \emph{independent} if no edge has both endpoints in $X$. The cardinality of a maximum independent set in $G$, written $\alpha(G)$, is the \emph{independence number} of $G$. The \emph{core} of $G$, written $\core(G)$, is the intersection of all maximum independent sets in $G$. 

The graph parameters $\alpha(G)$ and $\mu(G)$ are in general negatively correlated (adding edges doesn't increase the independence number and doesn't decrease the matching number) but incomparable as can be seen by the following observations.
\begin{ob}
If $G = E_n$, the $n$-vertex empty graph, then $0 = \mu(G) < \alpha(G) = n$.
\end{ob}

\begin{ob}
If $G = K_n$, the $n$-vertex complete graph, then $1 = \alpha(G) < \mu(G) = \lfloor \frac{n}{2} \rfloor$.
\end{ob}

However from the point of view of ``almost all graphs'', random graph theory provides us with high probability the bound $c_1\log(\mu(G)) \le \alpha(G)\le c_2\log(\mu(G))$~\cite{Bo2001, CoEf2015, FrPi2004}. Thus, with high probability in a random graph, $\mu(G)$ is much higher than $\alpha(G)$. 

Despite the above observations and examples, there exists many relationships between $\alpha(G)$ and $\mu(G)$. The following inequality is one of the most well known examples.
\begin{equation}\label{boros}
n(G) - 2\mu(G) \le n(G) - 2\mu^{*}(G)  \le \alpha(G) \le n- \mu(G)
\end{equation}
Here $\mu^{*}(G)$ denotes the cardinality of a minimum maximal matching in $G$. Graphs that satisfy the righthand side of (\ref{boros}) with equality are called \emph{König--Egerváry}, and have been extensively studied; see for example~\cite{BoDoDuFaGrSa13, LaPe2011, LeMa03, LeMa13}. Boros et al.~\cite{BoGoLe02} proved  $\alpha(G) \le \mu(G) + |\core(G)| - 1$ whenever $G$ is a graph with $\alpha(G) > \mu(G)$. Recently Levit et al.~\cite{LeMa19} proved a similar result, namely $\alpha(G) \le \mu(G) + |\core(G)| - |N_G(\core(G))|$ whenever $G$ is a graph with a matching from $N_G(\core(G))$ into $\core(G)$. Intersecting maximum independent sets were also studied by Deniz et al.~\cite{DeLeMa19}, who showed $\alpha(G) \le \mu(G)$, provided $G$ contains two disjoint maximum independent sets. Levit et al.~\cite{LeMa19} also showed that $\alpha(G) \le \mu(G)$, under the condition that $G$ contains a unique odd cycle. 

So the first source of motivation in our paper is to try and obtain a deeper understanding of these kinds of inequalities relating $\alpha(G)$ and $\mu(G)$ via the cardnality of the core of $G$. Another source of motivation comes from the following example: If $G$ is the bipartite graph $K_{\delta, n-\delta}$, where $n \ge \delta \ge 1$, then $\delta(G) = \delta$ and $\Delta(G) = n-\delta$. Clearly $n-\delta = \Delta(G) = \alpha(G)$ and $\delta = \mu(G)$, and so, $\alpha(G) = \frac{\Delta(G)}{\delta(G)}\mu(G)$. Thus, a natural question arises, namely, is this the best possible upper-bound on $\alpha(G)$ in terms of the parameters $\mu(G)$, $\delta(G)$, and $\Delta(G)$?

Our main two theorems supply answers to the problems and motivation mentioned aboce. These two theorems are shown below. 
\begin{thm}
\label{thm:main_X}
If $G$ is a graph and $X$ is any intersection of maximum independent sets, then 
\[
\alpha(G)  \le \mu(G) + |X| - \mu(G[N_G[X]]),
\]
and this bound is sharp.
\end{thm}

\begin{thm}\label{thm:main_delta}
If $G$ is a graph, then
\[
\delta(G)\alpha(G) \le \Delta(G)\mu(G),
\] 
and this bound is sharp.
\end{thm}

As can easily be seen, these two results generalize and in many cases improve on many of the known relationships between $\alpha(G)$ and $\mu(G)$. The remainder of our paper is structured as follows. In Section~\ref{S:intersecting} we prove Theorem~\ref{thm:main_X}. In Section~\ref{S:delta} we prove Theorem~\ref{thm:main_delta}. Finally in Section~\ref{S:problems}, we give concluding remarks, suggestions for future work, and a new conjecture. 

For notation and terminology not found here, we refer the reader to West~\cite{West}. We will also make use of the standard notation $[k] = \{1, \dots, k\}$.

\section{Proof of Theorem 1}\label{S:intersecting}
In this section we prove Theorem 1. Before doing so, we will need the following two lemma's. 
\begin{lem}\label{l:main_one}
If $A$ is an independent set and $X$ is a maximum independent set, then there is a matching from $A - (A\cap X)$ to $X - A$ that saturates each vertex in $A - (A\cap X)$\footnote{We acknowledge that $A -(A\cap X)$ is equivalent to $A - X$. However, we use $A -(A\cap X)$ in place of $A - X$ because this view of the set difference becomes useful in subsequent proofs.} . 
\end{lem}
\proof Since $A - (A \cap X)$ and $X - A$ are independent sets, respectively, we first observe that $[A - (A \cap X)] \cup (X - A)$ induces a bipartite subgraph in $G$ with bipartitions $A - (A \cap X)$ and $X - A$. Let $H$ denote this induced subgraph and choose $S\subseteq A - (A \cap X)$ arbitrarily. If $|S| > |N_H(S)|$, then $(X\cup S) - N_H(S)$ forms an independent set in $G$ with cardinality larger than $X$. This is a contradiction since $X$ is a maximum independent set in $G$, and so, it must be the case that $|S| \le |N_H(S)|$. Since $S$ was chosen arbitrarily, Hall's Theorem implies there is a matching from $A - (A \cap X)$ to $X - A$ saturating $A - (A \cap X)$ proving the lemma. \qed 

Using Lemma~\ref{l:main_one}, we next prove a technical lemma that bounds the difference between the size of independent sets and the matchings numbers of their closed neighborhoods. 
\begin{lem}\label{l:main_two}
If $A$ is an independent set and $X$ is any intersection of maximum independent sets with $X\subseteq A$, then 
\[
|A| -  \mu(G[N_G[A]]) \le  |X| - \mu(G[N_G[X]]).
\]
\end{lem}
\proof Let $X = X_1\cap \dots \cap X_k$, where $X_i$ is a maximum independent set in $G$ for each $i\in [k]$, and let $A\subseteq V(G)$ be any independent set satisfying $X\subseteq A$. For notational convenience, let $X_0 = A$ and $A_r = \bigcap_{i=0}^{r-1}X_{i}$. Note $A_r$ is an independent set for all $r \in [k]$. By Lemma~\ref{l:main_one} there is a matching from $X_r - A_r$ to $A_r -(A_r \cap X_r)$ that saturates every vertex contained in $A_r -(A_r \cap X_r)$. Let $M_r$ denote one such matching for each $r \in [k]$. Note that edges in $M_j$ and $M_i$ will not share endpoints for any $i \neq j$ and $i, j \in [k]$. Furthermore, each edge in $M_r$ contains at least one endpoint in $A$, again for each $r\in [k]$. Thus, $M = M_1, \dots,  M_k$ is a matching in the induced subgraph $G[N_G[A]]$. 

Thus far we have only saturated vertices in $A_{k} = X_0 \cap \dots \cap X_{k-1}$. Let $Q$ be a maximum matching in $G[N_G[X]]$.
Next observe that $Q$ is edge independent from the matching $M_1 \cup M_2 \cup \dots \cup M_{k-1}$. Thus, $M = M_1 \cup M_2 \cup \dots \cup M_{k} \cup Q$ is a matching in $G[N_G[A]]$. This implies the following inequality
\vskip 0.05 cm
\[
\begin{array}{lcl}
\mu(G[N_G[A]]) & \ge & \displaystyle{|M|} \1 \\ 
& = & \displaystyle{\big|\bigcup_{i=1}^{k}M_i\big| + |Q|} \1 \\
& = & \displaystyle{\sum_{i=1}^{k}\Big(|A_i - (A_i \cap X_i)|\Big) + \mu(G[N_G[X]])}\1 \\
& = & \displaystyle{\sum_{i=1}^{k}\Big(|A_i| - |A_{i+1}|\Big) + \mu(G[N_G[X]])}\1 \\
\\
& = & \displaystyle{|A| - |X| + \mu(G[N_G[X]])},
\end{array}
\] 
\vskip 0.10 cm
\noindent where all the terms in the summation, except the first and last, cancel out because the summation in the inequality is a telescoping series. Rearranging the above inequality, we finish the proof of our lemma. \qed

With Lemma~\ref{l:main_two} we next prove Theorem 1. Recall its statement. 
\begin{duplicate}
If $G$ is a graph and $X$ is any intersection of maximum independent sets, then 
\[
\alpha(G)  \le \mu(G) + |X| - \mu(G[N_G[X]]),
\]
and this bound is sharp.
\end{duplicate}
\proof Let $X$ be an intersection on maximum independent sets, one of which is the set $A$. By Lemma~\ref{l:main_two}, we have 
\[
|A| -  \mu(G[N_G[A]]) \le  |X| - \mu(G[N_G[X]]).
\]

Since $N_G[A] = V(G)$ implies $\mu(G[N_G[A]]) = \mu(G)$, and since $|A| = \alpha(G)$, we obtain
\[
\alpha(G) -  \mu(G) \le  |X| - \mu(G[N_G[X]]).
\]
Rearranging the above inequality proves the inequality posed in the theorem. To see this inequality is sharp, see Example~\ref{e:main}. \qed 

\begin{remark} Let $Q$ be a set of maximum independent sets in a graph $G$, where $|Q| \ge 3$. In light of Theorem 1, it is natural to ask what number of elements in $Q$ together form the optimal intersection with respect to the upper bound on $\alpha(G)$ given by the theorem? The answer comes from Lemma 4,  $X$ is the intersection of all elements in $Q$ and $A$ is the intersection of two elements in $Q$, then $|A| -  \mu(G[N_G[A]]) \le |X| - \mu(G[N_G[X]])$. Rearranging, we obtain
\[
\alpha(G) \le \mu(G) + |A| - \mu(G[N_G[A]]) \le \mu(G) + |X| - \mu(G[N_G[X]]).
\]
Thus, every collection of three or more elements in $Q$ has a pair that yields a better bound on $\alpha(G)$.
\end{remark}

As a consequence of Remark 1, we obtain the following corollary. 
\begin{cor}
If $G$ is a graph with no unique maximum independent set, and $Q$ is the set of all maximum independent sets in $G$, then
\[
\alpha(G) \le \mu(G) + \min\{|A\cap B| - \mu(G[N_G[A\cap B]]): \text{ } A, B \in Q\}.
\]
\end{cor}

The following example gives an infinite family of graphs satisfying Theorem 1 with equality. Moreover, it also shows a family of graphs where any intersection of maximum independent sets will satisfy Theorem 1 with equality. 
\begin{example}\label{e:main}
With this example, we establish the inequality of Theorem 1 being sharp, and in doing so, we also show the existence of graphs where any choice of intersecting maximum independent sets satisfies the inequality with equality. Let $p$, $q$, and $r$ be non-negative integers with $p+r \ge 2$. Let $G(p,q,r)$ be the graph obtained by attaching each vertex of $G_1$ (the complete graph $K_p$ with a pendant attached to each vertex) to each vertex of $G_2 = K_q$, and then attaching each vertex of the empty graph with order $r$, denoted $G_3$, to every vertex of $G_2$. For the graph $G(p, q, r)$, observe:
\begin{itemize}
\item[A.] $\core(G(p, q, r)) = V(G_3)$.

\item[B.] $\alpha(G) = p + r$.

\item[C.] 
$
\mu(G) = \left\{
        \begin{array}{ll}
            p+q, & \quad \text{if } r \ge q. \\
            p + \frac{r+q}{2}, & \quad \text{if } r < q.
        \end{array}
    \right.
$

\item[D.]  If $X$ is any intersection maximum independent sets in $G(p, q, r)$, then 
\[
|X| - \mu(G[N_G[X]]) = \left\{
        \begin{array}{ll}
            r - q, & \quad \text{if }r \ge q. \\
            \frac{r - q}{2}, & \quad \text{if } r < q.
        \end{array}
    \right.
\]

\end{itemize}
With the above equations, if $X$ is any intersection of maximum independent sets in $G(p, q, r)$, then
\[
\alpha(G(p, q, r)) = \mu(G(p, q, r)) + |X| - \mu(G[N_{G(p, q, r)}[X]]).\footnote{Note:  $|X| - \mu(G[N_{G(p, q, r)}[X]]) < 0$ whenever $r < q$.}
\]
\end{example}

The infinite family of graphs given in Example~\ref{e:main} provide examples where any choice of intersecting maximum independent sets will satisfy Theorem 1 with equality. The graph presented in Figure 1 provides an example where no choice of intersecting maximum independent sets will satisfy Theorem 1 with equality. 
\begin{figure}[htb]
\begin{center}
\vskip 0.35 cm
\begin{tikzpicture}[scale=0.8,style=thick,x=1cm,y=1cm, =>stealth]
\def\vr{2.5pt} 
\path (-1, -2) coordinate (w1);
\path (-1, 0) coordinate (x1);
\path (1, -1) coordinate (z1);

\draw (z1)--(x1);
\draw(x1)--(w1);
\draw (z1)--(w1);

\draw (w1) [fill=black]circle (\vr);
\draw (x1) [fill=black]circle (\vr);
\draw (z1) [fill=black]circle (\vr);


\path (5, -2) coordinate (w);
\path (5, 0) coordinate (x);
\path (3, -1) coordinate (z);

\draw (z)--(x);
\draw(x)--(w);
\draw (z)--(w);

\draw(x) -- (x1);
\draw(w)--(w1);
\draw(z)--(z1);

\draw (w) [fill=black]circle (\vr);
\draw (x) [fill=black]circle (\vr);
\draw (z) [fill=black]circle (\vr);

\path (2, -2) coordinate (v1);
\path (2, 0) coordinate (v2);
\path (2, -1) coordinate (v3);

\draw (v1) [fill=black]circle (\vr);
\draw (v2) [fill=black]circle (\vr);
\draw (v3) [fill=black]circle (\vr);


\end{tikzpicture}
\end{center}
\vskip -0.5 cm
\caption{A graph $G$ where no choice of intersecting maximum independent sets satisfies Theorem 1 with equality.} 
\end{figure}
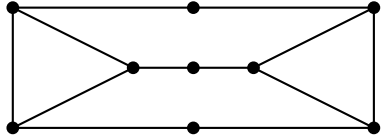

One interesting application of Theorem 1 can be seen by considering well-covered graphs, a heavily studied notion; see for example~\cite{FiHaNo1993, Ha1999, Pl1993}.  A graph is \emph{well-covered} if all maximal independent sets are also maximum. Observe that if $G$ is an isolate-free and well-covered graph, then for every vertex $v\in V(G)$ there is a neighbor of $v$, say $w$, so that $v$ cannot appear in any maximum independent set containing $w$. Since we may greedily construct a maximal independent set (which is also a  maximum independent set in well-covered graphs), starting from either $v$ or $w$, it follows that the intersection of all maximum independent sets in $G$ is necessarily empty. Therefore, taking $X = \core(G) = \emptyset$ in Theorem 1 implies the following corollary. 
\begin{cor}
If $G$ is an isolate-free and well-covered graph, then
\[
\alpha(G) \le \mu(G).
\]
\end{cor}

Theorem 1 also generalizes and improves several known results. For example, recall $\alpha(G) \le \mu(G)$, whenever $G$ contains two disjoint maximum independent sets (Deniz et al.~\cite{DeLeMa19}). Since Theorem 1 implies $\alpha(G) \le \mu(G)$ whenever \emph{any} collection of maximum independent sets has an empty intersection, their result is a corollary of Theorem 1. Another example comes from considering the bound $\alpha(G) \le \mu(G) + |\core(G)| - 1$, whenever $\alpha(G) > \mu(G)$ (Boros et al.~\cite{BoGoLe02}). Taking $X = \core(G)$ in the statement of Theorem 1, observe that if $\alpha(G) > \mu(G)$ and $\mu(G[N_G[\core(G)]]) > 1$, then Theorem 1 improves upon this result. In particular, we make note of the following corollary. 
\begin{cor}
If $G$ is a graph, then
\[
\alpha(G) \le \mu(G) + |\core(G)| - \mu(G[N_G[\core(G)]]).
\]
\end{cor}

\section{Proof of Theorem 2}\label{S:delta}
In this section we prove Theorem 2. Before doing so we first prove a theorem and recall a lemma. The following result was motivated by a conjecture of the automated conjecturing program TxGRAFFITI, which in turn was motivated by GRAFFITI of Fajtlowicz~\cite{Graffiti}, and later GRAFITTI.pc of DeLaVeña\cite{Graffiti.pc}. The program TxGRAFFITI was written by the second author, and generates possible graph inequalities on simple connected graphs. When asked to conjecture on the independence number, the program conjectured $\alpha(G) \le \mu(G)$ for all $3$-regular and connected graphs. The following theorem confirms and generalizes this conjecture. 
\begin{thm}\label{alpha_mu_thm}
If $G$ is a $r$-regular graph with $r >0$, then 
\[
\alpha(G) \le \mu(G).
\]
\end{thm}
\proof Let $G$ be an $r$-regular graph with $r>0$, $X\subseteq V(G)$ be a maximum independent set , and $Y = V(G) - X$. By removing edges from $G$ with both endpoints in $Y$, we next form a bipartite graph $H$ with partite sets $X$ and $Y$. Since those edges removed from $G$ in order to form $H$ were only edges with both endpoints in $Y$, any vertex chosen in $X$ will have the same open neighborhood in $H$ as it does in $G$. It follows that since $G$ is $r$-regular and since $X$ is an independent set, any vertex in $X$ will have exactly $r$ neighbors in $Y$, both in $G$ and in $H$.

Let $S\subseteq X$ be chosen arbitrarily, and let $e(S, N_{H}(S))$ denote the number of edges from $S$ to $N_{H}(S)$. Since each vertex in $S$ has exactly $r$ neighbors in $Y$, we observe that $e(S, N_{H}(S)) = r|S|$. However, since each vertex in $N_{H}(S)$ has at most $r$ neighbors in $X$, we also have $e(S, N_{H}(S)) \le r|N_{H}(S)|$. It follows that $r|S| \le r|N_{H}(S)|$, and so, $|S| \le |N_{H}(S)|$. By Hall's Theorem~\cite{Ha35}, there exists a matching $M$ that can match $X$ to a subset of $Y$. Since $X$ is a maximum independent set and since $M$ is also a matching in $G$, we conclude $\alpha(G) = |M| \le \mu(G)$, proving the theorem. \qed  

A $k$-edge-coloring of $G$ is an assignment of $k$ colors to the edges of $G$ so that no two edges with the same color share an endpoint. The minimum integer $k$ so that $G$ has a $k$-edge-coloring is the \emph{edge chromatic number} of $G$, written $\chi'(G)$. By Vizing's Theorem, $\Delta(G) \le \chi'(G) \le \Delta(G) + 1$ for all graphs. Graphs satisfying $\Delta(G) = \chi'(G)$ are \emph{class 1}, whereas graphs satisfying $\chi'(G) = \Delta(G) + 1$ are \emph{class 2}. Let $G_\Delta$ denote the subgraph induced by the set of maximum degree vertices in $G$. With these definitions, we next recall a result due to Fournier~\cite{fournier}.  
\begin{lem}\label{l:fournier}
If $G$ is class 2, then $G_\Delta$ contains at least one cycle. 
\end{lem}

As a consequence of Lemma~\ref{l:fournier}, all class 2 graphs satisfy $|G_\Delta| \ge 3$ and $E(G_\Delta) \neq \emptyset$. With this observation, we are now ready to prove Theorem 2. Recall its statement. 
\begin{duplicate}
If $G$ is a graph, then
\[
\delta(G)\alpha(G) \le \Delta(G)\mu(G),
\] 
and this inequality is sharp. 
\end{duplicate}
\proof Clearly, if $\delta(G) = 0$ we are done. So we will assume $\delta(G) > 0$. Proceeding by way of contradiction, suppose the theorem is false. Among all counter-examples, let $G$ be one with a minimum number of edges. By Theorem~\ref{alpha_mu_thm}, any $r$-regular graph with $r>0$ will satisfy the theorem, and so, the graph $G$ must satisfy $\delta(G) < \Delta(G)$. Before proceeding, we remind the reader that all graphs are either class 1 or class 2. 

If $G$ is a class 1, then $\chi'(G) = \Delta(G)$. Since each color class in a $\chi'(G)$-edge coloring forms a matching in $G$, and since every edge in $G$ belongs to exactly one color class, it is clear that $m(G) \le \chi'(G)\mu(G)$. Moreover, each vertex in any maximum independent set will have at least $\delta(G)$ edges incident with it, implying $\delta(G)\alpha(G) \le m(G)$. It follows that $\delta(G) \alpha(G) \le \Delta(G)\mu(G)$, which is impossible, because $G$ is a counter-example. Thus, $G$ is not class 1. 

If $G$ is class 2, then Lemma~\ref{l:fournier} implies $G_\Delta$ has a non-empty edge set. Let $vw$ be one such edge and let $H = G - vw$. Clearly, $\alpha(G) \le \alpha(H)$ and $\mu(H)\le \mu(G)$. Since $\delta(G) < d_G(v) = d_G(w) = \Delta(G)$, it follows that $\delta(H) = \delta(G)$. Since $G_\Delta$ contains a cycle, it has at least 3 vertices, and so, $\Delta(H) = \Delta(G)$. Finally, $G$ being a minimum counterexample implies $\delta(H)\alpha(H) \le \Delta(H)\mu(H)$. It follows that $\delta(G)\alpha(G) \le \Delta(G)\mu(G)$, which is again impossible, since $G$ is a counter-example. Thus, $G$ is not class 2. Since $G$ is neither class 1 nor class 2, we contradict the existence of $G$.  

To see this bound is sharp, first consider $(\delta ,\Delta)$-bipartite graphs. Namely, the graph $G$ with $V(G) = A\cup B$ where  $A$ and $B$ are independent sets and all degrees in $A$ equal $\Delta(G)$ and all degrees in $B$ equal  $\delta(G)$. Thus showing the bound sharp for class 1 graphs. \qed

We next consider applications of Theorem 2. More specifically, if $G$ is a graph with $\delta(G) \ge 1$, then Theorem 2 implies 
\begin{equation}
\alpha(G) \le \frac{\Delta(G)}{\delta(G)}\mu(G).
\end{equation}
This bound is interesting, as the righthand side of (2) is computable in polynomial time. Moreover, (2) can also improve on known computationally efficient upper bounds for $\alpha(G)$ in some classes of graphs. For example, the \emph{annihilation number} of $G$, written $a(G)$, is a degree sequence invariant for which $\alpha(G) \le a(G)$~\cite{Pe04, Pe09}. This bound improves on many known bounds, for example $\alpha(G) \le a(G) \le n(G) - \frac{m(G)}{\Delta(G)}$ (see~\cite{Pe04}). However, $a(G) \ge \frac{n(G)}{2}$ for all graphs. Thus, for $r$-regular graphs with $r >0$, (2) gives the improvements
\[
\alpha(G) \le \mu(G) \le \frac{n(G)}{2} \le a(G)\le n(G) - \frac{m(G)}{\Delta(G)},
\]
and
\[
\alpha(G) \le \mu(G) \le n(G) - \mu(G).
\]
Further observe that for sufficiently large $r$-regular graphs with $r>0$, (2) can give dramatic improvements on the minimum degree bound $\alpha(G) \le n(G) - \delta(G)$.

\section{Concluding Remarks}\label{S:problems}
In this paper we have proven two theorems relating $\alpha(G)$ and $\mu(G)$. These two theorems imply a myriad of interesting corollaries bounding $\alpha(G)$ from above; some of which we summarize in the following theorem.
\begin{thm}\label{thm:last}
Let $G$ be a graph and $X\subseteq V(G)$ be any intersection of maximum independent sets.
\begin{itemize}
\item[1.] $\alpha(G)  \le \mu(G) + |X| - \mu(G[N_G[X]])$.

\item[2.] $\alpha(G) \le \mu(G) + \core(G) - \mu(G[N_G[\core(G)]])$.

\item[3.] If $\core(G) = \emptyset$, then $\alpha(G) \le \mu(G)$.

\item[4.] If $X$ is isolate-free and well-covered, then $\alpha(G) \le \mu(G)$.

\item[5.] If $\delta(G) \ge 1$, then $\alpha(G) \le \frac{\Delta(G)}{\delta(G)}\mu(G)$.

\item[6.] If $G$ is $r$-regular with $r>0$, then $\alpha(G) \le \mu(G)$. 

\end{itemize}

\end{thm}

As mentioned before, many of the cases contained in Theorem~\ref{thm:last} yield improvements on known upper bounds for $\alpha(G)$, most notably being the case of Theorem~\ref{thm:last}.6. Observing this, we believe the following problem merits further inspection. 
\begin{prob}
Characterize $\alpha(G) = \mu(G)$ whenever $G$ is 3-regular. 
\end{prob}
More generally, we also suggest the following problem. 
\begin{prob}
Characterize all graphs $G$ for which $\delta(G)\alpha(G) = \Delta(G)\mu(G)$. 
\end{prob}

Next we remark on the cardinality of minimum maximal matchings in $G$, written $\mu^*(G)$. Recall, 
\[
n(G) - 2\mu^*(G) \le \alpha(G) \le n(G) - \mu(G),
\]
for any graph $G$. Thus, by Theorem~\ref{thm:last} we obtain
\begin{equation}\label{eq:mu}
n(G) - 2\mu^*(G) \le \mu(G),
\end{equation}
whenever $G$ is $r$-regular with $r>0$, or isolate-free and well-covered, or has an empty core. Rearranging~(\ref{eq:mu}), we obtain the inequality
\[
\frac{\alpha(G)}{2} \le  \frac{n(G) - \mu(G)}{2} \le  \mu^*(G),
\]
for all graphs satisfying one or more of the above mentioned properties. Hence, $\alpha(G) \le 2\mu^*(G)$ for these families of graphs. We suggest that future work include studying relationships between independent sets and $\mu^*(G)$. More specifically, we suggest considering the following conjecture of TxGRAFFITI.
\begin{conj}\label{i_mu_conj}
If $G$ is an $r$-regular graph with $r >0$, then 
\[
i(G) \le \mu^*(G),
\]
where $i(G)$ denotes the minimum cardinality among all maximal independent sets\footnote{This graph invariant is known as the \emph{independent domination number}, and has also been heavily studied in the literature; see for example the excellent survey by Goddard and Henning~\cite{GoHe2013}.} in $G$.
\end{conj}
If Conjecture~\ref{i_mu_conj} is true, we believe the following question merits further investigation.
\begin{quest}
Is it true that $\delta(G)i(G) \le \Delta(G)\mu^*(G)$ for all graphs?
\end{quest}

Finally, we would like to acknowledge and thank Craig Larson for his early conversations on the some of the conjectures of TxGRAFFITI presented in this paper.

\end{document}